\documentclass[11pt]{amsart}
\usepackage{epsfig}
\newcommand{\bT}{\bar{T}}

\newtheorem{thm}{Theorem}[section]
\newtheorem{cor}[thm]{Corollary}
\newtheorem{lemma}[thm]{Lemma}
\newtheorem{prop}[thm]{Proposition}

\newtheorem{q}[thm]{Question}
\theoremstyle{definition}
\newtheorem{defn}[thm]{Definition}

\begin{document}

\title{Graph Powers and $k$-Ordered Hamiltonicity}

\author{Denis Chebikin}
\thanks{Please send all correspondence to: 1404 Yorkshire Ln,
  Shakopee, MN 55379}
\address{Department of Mathematics, Massachusetts Institute of
  Technology, Cambridge, MA 02139}
\email{chebikin@mit.edu}


\begin{abstract}
It is known that if $G$ is a connected simple graph, then $G^3$ is
Hamiltonian (in fact, Hamilton-connected). A simple graph is
$k$-ordered Hamiltonian if for any sequence $v_1$, $v_2$, \dots, $v_k$
of $k$ vertices there is a Hamiltonian cycle containing these vertices
in the given order. 
In this paper, we prove that $G^{\lfloor 3k/2 \rfloor +1}$ is
$k$-ordered Hamiltonian for a connected graph $G$ on at least $k$ vertices.
We further show that if $G$ is connected, then $G^4$ is
4-ordered Hamiltonian and that if $G$ is Hamiltonian, then $G^3$ is
5-ordered Hamiltonian. We also give bounds on the smallest power $p_k$
such that $G^{p_k}$ is $k$-ordered Hamiltonian for $G=P_n$ and $G=C_n$.
\end{abstract}

\maketitle

\begin{center}
\textit{Dedicated to Pavlo Pylyavskyy on the occasion of his 21st birthday.}
\end{center}

\section{Introduction}\label{intro}

The concept of $k$-ordered Hamiltonian graphs has been recently
introduced by Ng and Schultz \cite{Ng}. A simple graph $G$ is
$k$-\textit{ordered} (resp. $k$-\textit{ordered Hamiltonian}) if for
any sequence $v_1$, $v_2$, \dots, $v_k$ of $k$ vertices of $G$ there
is a cycle (resp. a Hamiltonian cycle) in $G$ containing these
vertices in the given order. Note that being $3$-ordered Hamiltonian
is equivalent to being Hamiltonian.

A natural direction of research related to this new Hamiltonian
property is to generalize existing results implying graph
Hamiltonicity and obtain results implying $k$-ordered
Hamiltonicity. In \cite{Ng} Ng and Schultz generalize classical
theorems of Dirac and Ore and give minimum vertex degree conditions
that guarantee a graph is $k$-ordered Hamiltonian. These
conditions were improved by Faudree \cite{Faudree}. Another series of
results appearing in \cite{Faudree2} describes 
various forbidden subgraphs that force a
graph to be $k$-ordered or $k$-ordered Hamiltonian. There are many
open questions about whether these properties can be ensured by
sufficient connectivity in a graph (see \cite{Faudree}).

In this paper we extend a well-known result on
Hamiltonicity of the third power of a simple graph, which is defined below.

\begin{defn}\label{g_power}
Let $G$ be a simple graph with vertex set $V$ and edge set $E$. 
The $n$-\textit{th power} of $G$, denoted $G^n$, is the simple graph
with the same vertex set $V$ and with the edge set
$$
E(G^n) = \{(v,w)\ |\ d_G(v,w) \leq n\}.
$$
\end{defn}

Our goal is to explore the $k$-ordered
Hamiltonicity of graph powers. 
In Section \ref{ppaths} we give bounds on the smallest power of the
path $P_n$ that is $k$-ordered Hamiltonian. In Section \ref{main} we
prove the main theorem, which states that $G^{\lfloor 3k/2 \rfloor
  +1}$ is $k$-ordered Hamiltonian for a connected graph $G$ on at
least $k$ vertices.
In
Section \ref{four} we prove the 4-ordered Hamiltonicity of the fourth
power of a connected graph.  In Section \ref{cycles} we
discuss the $k$-ordered Hamiltonicity of powers of the cycle $C_n$.
In Section \ref{conclusion} we pose several open questions related to
our results.

\section{Preliminaries}\label{prelim}

All graphs considered in this paper are finite simple graphs. The distance
between vertices $v$, $w$ of a graph $G$ is denoted by $d_G(v,w)$.
The number of vertices of $G$ is denoted by $|G|$. If $P$ is a path
with endpoints $a$ and $b$, then $aPb$ denotes the path $P$ starting
at $a$ and ending at $b$.

A graph $G$ is \textit{Hamilton-connected} if for any pair $v,w$ of
vertices of $G$ there is a path in $G$ between $v$ and $w$ containing
all vertices of $G$. Such path is called a \textit{Hamiltonian path}.
The following theorem is often included as an exercise 
in graph theory textbooks.

\begin{thm}\label{g3c}
If $G$ is a connected graph on 2 or more vertices, then $G^3$ is
Hamilton-connected.
\end{thm}

\begin{proof}
Let $T$ be a spanning tree of $G$. Notice that it suffices to show
that $T^3$ is Hamilton-connected since $T^3$ is a subgraph of $G^3$
with the same vertex set.

We will show by induction that if $T$ is a tree, then $T^3$ is
Hamilton-connected. This is obvious if $T$ has only 2
vertices. 

Suppose that the assertion is true for trees with fewer than
$|T|$ vertices. Let $v_1$ and $v_2$ be distinct vertices of
$T$. Choose an edge $e = (w_1, w_2)$ of $T$ such that $T-e$ consists
of two connected components $T_1$ and $T_2$ satisfying $v_1,w_1\in
T_1$ and $v_2,w_2\in T_2$. For $i=1,2$ let $u_i = w_i$ if $w_i\neq
v_i$ or if $|T_i| = 1$,
otherwise let $u_i$ be a vertex of $T_i$ such
that $d_T(u_i,w_i) = 1$. Let $P_i$ be a Hamiltonian path in $T_i^3$
between
$v_i$ and
$u_i$ (if $u_i = v_i = w_i$, then
$P_i$ is the ``empty'' path starting and ending at $v_i$). Notice that
$d_T(u_1, u_2) \leq 3$, so $v_1P_1u_1u_2P_2v_2$ is a Hamiltonian path
in $T^3$ between $v_1$ and $v_2$.
\end{proof}

Since Hamilton-connectedness implies Hamiltonicity, we have the following
immediate corollary:

\begin{cor}\label{g3}
If $G$ is a connected graph on 3 or more vertices, then $G^3$ is Hamiltonian.
\end{cor}

The proof of Theorem \ref{g3c} has been included for two
reasons. First, it uses the fact that it suffices to prove the
statement only for trees instead of arbitrary graphs. We will use the
same idea in the proof of Theorem \ref{g4}. Second, it is based on a
simple induction argument, which cannot be applied if we need to keep
track of the order in which certain vertices are contained in the cycle.

The following theorem about 2-connected graphs was proved by Fleischner
\cite{Fleischner} in 1974. A simpler proof can be found in \cite{Diestel}.

\begin{thm}\label{g2}
If $G$ is a 2-connected graph on 3 or more vertices, then $G^2$ is Hamiltonian.
\end{thm}

Another result relevant to the discussion of this paper was proved by
Lou, Xu, and Yao \cite{Lou}.

\begin{thm}\label{2edges}
Let $G$ be a connected graph on 3 or more vertices. For any two edges
$e_1$ and $e_2$ of $G$, there is a Hamiltonian cycle in $G^4$ containing
$e_1$ and $e_2$. 
\end{thm}

For $k\geq 3$ and a graph $G$, let $p_k(G)$ be the
smallest integer
$p$ such that $G^p$ is $k$-ordered Hamiltonian. In this paper we give bounds
on $p_k(G)$ for an arbitrary connected graph $G$ and in the special cases of
$G$ being a path or a cycle.

\section{Powers of paths}\label{ppaths}

Let $P_n$ be the path on $n$ vertices. In this section we present
bounds on $p_k(P_n)$.

\begin{thm}\label{pathlowerbound}
For $k\geq 3$ and $n\geq 2k-1$, the
inequality $p_k(P_n)
\geq \left\lfloor{3k\over 2}\right\rfloor - 2$ holds.
\end{thm}

\begin{proof}
First, consider the case of even $k$, and let $k=2m$. 
We show that $(P_n)^{3m-3}$ is not $2m$-ordered.
Let
$v_1$, $v_3$, \dots, $v_{2m-1}$, $u_1$, $u_2$, \dots,
$u_{2m-1}$,
$v_2$, $v_4$, \dots, $v_{2m}$ be consecutive vertices
of the path
$P_n$. 
Suppose that $C$ is a cycle in $(P_n)^{3m-3}$
containing the vertices
$v_1$, $v_2$, \dots, $v_{2m}$ in order.
For $1\leq i \leq 2m$, let $R_i$ denote the part of
$C$ between, and including, $v_i$ and $v_{i+1}$
(indices taken
modulo $2m$). 
Put $U = \{u_1,u_2,\dots,u_{2m-1}\}$. Since $|U|<2m$,
there is an
index $i$ such that $R_i$ does not contain a vertex of
$U$. In $P_n$, the
set $U$ is located between $v_i$ and
$v_{i+1}$,
so $R_i$ must contain an edge $e$
that ``skips over'' 
$U$.
The edge $e$ connects one of $v_1$, $v_3$, \dots,
$v_{2m-1}$ and
one of $v_2$, $v_4$, \dots, $v_{2m}$ because
$d_{P_n}(v_1,v_2) =
d_{P_n}(v_{2m-1},v_{2m}) = 3m-1 > 3m-3$. Also, the
edge $e$ connects
two non-consecutive vertices of $v_1$, $v_2$, \dots,
$v_{2m}$ since
the distance in $P_n$ between two consecutive vertices
of this
sequence is at least $3m-2$. Thus, the cycle $C$
cannot contain $e$,
which contradicts $R_i$ containing $e$.

We treat the case of odd $k$ in a similar way. Let
$k=2m+1$. We
show that $(P_n)^{3m-2}$ is not $(2m+1)$-ordered. Let
$v_1$, $v_3$,
\dots, $v_{2m+1}$, $u_1$, $u_2$, \dots, $u_{2m-1}$,
$v_2$, $v_4$,
\dots, $v_{2m}$ be consecutive vertices of the path
$P_n$. Suppose
that $C$ is a cycle in $(P_n)^{3m-2}$ containing the
vertices
$v_1$, $v_2$, \dots, $v_{2m+1}$ in order. As in the
previous case,
we find an edge $e$ of $C$ connecting one of $v_1$,
$v_3$, \dots,
$v_{2m+1}$ and one of $v_2$, $v_4$, \dots, $v_{2m}$.
Again, this
edge connects two non-consecutive vertices of $v_1$,
$v_2$, \dots,
$v_{2m+1}$ since $d_{P_n}(v_i,v_{i+1}) \geq 3m-1$ for
$1\leq i \leq
2m$. We obtain a similar contradiction.

We conclude that $p_{2m}(P_n) \geq 3m-2$ and
$p_{2m+1}(P_n) \geq
3m-1$. The theorem follows.
\end{proof}

Next, we present a construction which shows that $P_n^{\lfloor 3k/2
  \rfloor - 1}$ is $k$-ordered Hamiltonian.

\begin{thm}\label{pathupperbound}
For $n\geq k\geq 3$, the
inequality
$p_k(P_n) \leq \left\lfloor{3k\over 2}\right\rfloor
-1$ holds.
\end{thm}

\begin{proof}
Let $t = \left\lfloor{3k\over 2}\right\rfloor-1$. We
show that
$(P_n)^t$ is $k$-ordered Hamiltonian. 
Label consecutive vertices of $P_n$ with consecutive
integers between
$1$ and $n$.
Let $v_1$, $v_2$, \dots, $v_k$
be a sequence of $k$ vertices of $P_n$. 
We view the $v_i$'s as both numbers and vertices.

First, we construct a cycle $C$ in $(P_n)^t$
containing the vertices
$v_1$, $v_2$, \dots, $v_k$ in order. We mark some of
the $v_i$'s using the following procedure:

\begin{enumerate}

\item\label{Step1}
Initially, put $V = \{v_1,\dots,v_k\} - \{v_\ell,
v_r\}$, where
$v_\ell$ and $v_r$ are the smallest and the largest
among the $v_i$'s,
respectively. Mark each $v_i\in V$ for which no other
element of
$V$ is congruent to $v_i$ modulo $t$.

\item\label{Step2}
Now, add $v_\ell$ to $V$, and mark $v_\ell$ if no
other element of $V$ is
congruent to $v_\ell$ modulo $t$.

\item\label{Step3}
Repeat step \ref{Step2} with $v_r$. (Note that if
$v_\ell \equiv v_r
\pmod t$ and $v_\ell$ is marked at step \ref{Step2},
then $v_r$ is not
marked.)

\end{enumerate}

Let $S(V)$ denote the set of residues modulo $t$
present
among elements of $V$. Also, let $u(V)$ be the number
of elements of
$V$ that are not marked. 

\begin{lemma}\label{justalemma}
After the above procedure, the inequality
$$
|S(V)| + u(V) \leq t
$$
holds.
\end{lemma}

\begin{proof}[Proof of Lemma \ref{justalemma}]
Consider the situation immediately after step \ref{Step1}, when $|V|
= k-2$.
Let $s$ be the number of residues modulo $t$ that
appear exactly
once among the elements of $V$, and let $m$ be the
number of residues
that appear more than once. Then $|S(V)| = s+m$ and
$u(V) = k-2-s$
since exactly $s$ elements are marked at step
\ref{Step1}. Thus
$$
|S(V)| + u(V) = k-2+m.
$$
Each of the $m$ ``repeated'' residues appears at least
twice, so
$m \leq \left\lfloor{k-2\over 2}\right\rfloor$.
Substituting this
relation into the above equation, we get
\begin{equation}\label{eq13}
|S(V)|+u(V) \leq \left\lfloor{3(k-2)\over 2}\right\rfloor.
\end{equation}

At step \ref{Step2}, either $v_\ell$ adds a new
residue to $V$, or
$v_\ell$ is not marked, but not both. Thus, the sum
$|S(V)|+u(V)$ is
increased by $1$. Similarly, at step \ref{Step3} the
sum $|S(V)|+u(V)$
is increased by $1$ again. Hence after the entire
procedure, we
have
$$
|S(V)| + u(V) \leq \left\lfloor{3(k-2)\over
2}\right\rfloor + 2 =
\left\lfloor{3k\over 2}\right\rfloor -1 = t,
$$
where $V$ is now the set of all the $v_i$'s.
\end{proof}

Let $v_{j_1}$,\dots, $v_{j_{u(V)}}$ be the vertices of
$V$ that are
not marked. Select $u(V)$ distinct
residues $r_1$, \dots, $r_{u(V)}$ modulo $t$ not
contained in $S(V)$.
Note that we can choose these residues because $|S(V)|
+ u(V) \leq t$. Set
$$
r(v_i) = \left\{\begin{array}{ll}
                                  v_i\bmod t, & \mbox{if $v_i$ is
				    marked,}\\
                                  r_{\ell}, & \mbox{if $v_i$ is not
				    marked and $v_i = v_{j_\ell}$}.
\end{array}\right.
$$
Note
that the
residues $r(v_1)$, \dots, $r(v_k)$ are all distinct
since at most
one instance of each residue modulo $t$ in $V$ is
marked.

For $1\leq i \leq k$, let $R_i$ be the path in
$(P_n)^t$ from $v_i$
to $v_{i+1}$ whose interior vertices are all the
vertices of $P_n$
between $v_i$ and $v_{i+1}$ congruent to $r(v_i)$
modulo $t$ (indices
taken modulo $k$). Let $C$ be the cycle obtained by
joining the paths
$R_1$, $R_2$, \dots, $R_k$.

Clearly, $C$ contains the vertices $v_1$, \dots, $v_k$
in order. We
need to check that the paths $R_i$ are interior
vertex-disjoint
and that no path $R_i$ contains a vertex of $V$ in its
interior.
The former property of the $R_i$'s follows from the
fact that the
values of $r(v_i)$ are all distinct. We check the
latter property.
Suppose that $R_i$ contains a vertex $v_j\in
V-\{v_\ell,v_r\}$ 
in its interior. Then
$$v_j \equiv r(v_i) \pmod t.$$ 
If $v_i$ is marked, then $r(v_i)=v_i$, so, by the
above congruence,
$v_i$ cannot be marked. If $v_i$ is not marked, 
then $r(v_i) \notin S(V)$,
which contradicts the above congruence.

Finally, we extend $C$ to a Hamiltonian cycle in
$(P_n)^t$. If $C$
does not contain all vertices of $P_n$, choose a
vertex $z\notin C$
adjacent in $P_n$ to a vertex $x\in C$. 
Let $w$ and $y$ be the two vertices
adjacent to $x$ in $C$. It is easy to check that
either
$d_{P_n}(z,w)$ or $d_{P_n}(z,y)$ is at most $t$.
Without loss of
generality, let $d_{P_n}(z,w)\leq t$. We can add $z$
to
$C$ by replacing the edge $(w,x)$ with the path $wzx$
to produce
a new cycle $C'$ in $(P_n)^t$. Clearly, $C'$ still
contains the
vertices $v_1$, \dots, $v_k$ in order. In this way we can
add all 
vertices of $P_n$ to the cycle.
\end{proof}

\section{Main theorem}\label{main}

We now prove our main result which gives an upper bound on $p_k(G)$
for a connected graph $G$ on at least $k$ vertices. We
begin by proving three technical lemmas.

\begin{lemma}\label{tl0}
Let $W$ be a tree on at least $2$ vertices, and let $w$ be a vertex of
$W$. Then $W^3-w$ has a Hamiltonian path whose endpoints $w_1$ and
$w_2$ satisfy $d_W(w,w_1)=1$ and $d_W(w,w_2)\leq 2$.
\end{lemma}

\begin{proof}
Let $W_1$, \dots, $W_m$ be the connected components of $W-w$. Let
$w^{(i)}_1$ be the vertex of $W_i$ adjacent to $w$ in $W$, and let
$w^{(i)}_2$ be a vertex of $W_i$ adjacent to $w^{(i)}_1$, or
equal to $w^{(i)}_1$ if $|W_i|=1$. Let $R^{(i)}$ be a Hamiltonian path
in $W_i^3$ starting at $w^{(i)}_1$ and ending at $w^{(i)}_2$; such a
path exists by Theorem \ref{g3c}. Note
that $d_W\bigl(w^{(i)}_2, w^{(i+1)}_1\bigl) \leq d_W\bigl(w^{(i)}_2,
w\bigl) + d_W\bigl(w, 
w^{(i+1)}_1\bigl) \leq 3$ for $1\leq i \leq m-1$, so the concatenation 
$R^{(1)}\dots R^{(m)}$ yields the desired path with $w_1 = w^{(1)}_1$
and $w_2 = w^{(m)}_2$.
\end{proof}

\begin{lemma}\label{tl1}
Let $k$ and $p\geq 3$ be positive integers, 
and let $G$ be a connected graph on at least $k$
vertices. Suppose that, for every sequence $v_1,\dots,v_k$ of
$k$ vertices of $G$, there exists a subtree $U\subseteq G$ and a cycle $C$ 
in $U^p$ satisfying the following
conditions:
\begin{enumerate}

\item[(i)] $C$ contains $v_1$, \dots, $v_k$ in order;

\item[(ii)] $C$ contains every leaf of $U$;

\item[(iii)]
if $x$ is a leaf of $U$, then $x$ is adjacent in $C$ to a
  vertex $y$ such that either $d_U(x,y)\leq p-2$, or $d_U(x,y)=p-1$ and
  $y$ is not a leaf of $U$.
\end{enumerate}
Then $G^p$ is $k$-ordered Hamiltonian.
\end{lemma}

\begin{proof}
We begin by extending $C$ to a Hamiltonian cycle in $U^p$. Let $y$ be
a vertex of $U$ such that $y\notin C$. By condition (ii), $y$ is not a
leaf. Therefore, $U-y$ has at least
two connected components, and some part 
of $C$
between a leaf in one component and a leaf in another component
does not contain $y$. This part
contains an edge $(x,z)$ such that $y$ lies on the unique path in $U$
between $x$ and $z$. Let $C'$ be the cycle obtained from $C$ by
replacing the edge $(x,z)$ with the path $xyz$. 
Clearly, $C'$ satisfies conditions (i) and (ii).
Since
$d_U(x,y)+d_U(y,z)=d_U(x,z)\leq p$, it follows that $d_U(x,y) \leq p-1$
and $d_U(z,y)\leq p-1$, hence $C'$ satisfies condition (iii). In this
way, we can add all remaining vertices of $U$ to obtain a Hamiltonian
cycle $\tilde{C}$ of $U^p$ satisfying conditions (i)---(iii) with $C$
replaced by $\tilde C$.

Let $T$ be a spanning tree of $G$ containing $U$.
Consider the graph $T-E(U)$ obtained by
removing the edges of $U$ from $T$. Let $U_1$,
\dots, $U_\ell$ be the connected components of $T-E(U)$ such that 
$|U_i|>1$. Let $u_i$ denote the unique vertex
of
$U\cap U_i$. Note that the vertices $u_1$, \dots, $u_\ell$ are
distinct.
Put $T_0 = U$ and $T_i = U \cup
U_1 \cup \dots \cup U_i$ for $1\leq i \leq \ell$, so that $T_\ell = T$. 
We construct a
sequence $\tilde C = C_0$, $C_1$, \dots, $C_\ell$, where $C_i$ is a
Hamiltonian cycle of $T_i^p$ satisfying the following conditions:
\begin{enumerate}
\item[(I)] $C_i$ contains the vertices $v_1$, \dots, $v_k$ in order;

\item[(II)] if $x$ is a common leaf of $U$ and $T_i$, then $x$ is
  adjacent in $C_i$ to a vertex $y$ such that either $d_T(x,y) \leq
  p-2$, or $d_T(x,y) = p-1$ and $y$ is not a leaf of $U$.
\end{enumerate}

The cycle $C_0$ satisfies the above conditions. We construct $C_{i+1}$
by inserting the vertices of $U_{i+1}-u_{i+1}$ between two consecutive
vertices of $C_i$. By Lemma \ref{tl0}, there is a Hamiltonian path
$w_1R_iw_2$
in
$U^3_{i+1}-u_{i+1}$ whose endpoints $w_1$ and $w_2$ satisfy
$d_T(u_{i+1},w_1)=1$ and $d_T(u_{i+1},w_2)\leq 2$.
We consider two cases.

\vskip5pt

Case 1: $u_{i+1}$ is not a leaf of $T_i$. Then $T_i-u_{i+1}$ has at least
two
connected components, and one of the parts of $C_i$ between a vertex
in one component and a vertex in another other component does not contain
$u_{i+1}$. This part contains an edge $(x,z)$ such that $u_{i+1}$
belongs to the unique path between $x$ and $z$ in $T_i$. Without loss
of generality, suppose that $d_T(x,u_{i+1}) \geq d_T(z,u_{i+1})$. Let
$C_{i+1}$ be the cycle obtained by replacing the edge $(x,z)$ with the path
$xw_1R_iw_2z$ in $C_i$. Then $C_{i+1}$ is a Hamiltonian cycle in
$T_{i+1}^p$ because $d_T(x,w_1) = d_T(x,u_{i+1}) + 1 \leq p$ and
$d_T(z,w_2) \leq d_T(z,u_{i+1})+d_T(u_{i+1},w_2)
 \leq \lfloor p/2 \rfloor + 2 \leq p$.

Since $C_{i+1}$ contains the vertices of $C_i$ in the same order as
$C_i$, it follows that $C_{i+1}$ satisfies condition (I). To show that
$C_{i+1}$ satisfies condition (II), we need to consider the case when
$x$ or $z$ is a common leaf of $U$ and $T_{i+1}$. If $x$ is a common
leaf of $U$ and $T_{i+1}$, then $x$ is also a leaf of $T_i$, so either
$d_T(x,z) \leq p-1$, or the vertex $y\neq z$ adjacent to
$x$ in $C_i$ has the property of condition (II). In the former case, we have
$d_T(x,w_1) = d_T(x,u_{i+1})+1 \leq p-1$, and $w_1$ is not a leaf of $U$. 
In the latter case, the
vertex $y$, which is adjacent to $x$ in $C_{i+1}$, has the required
property. The case when $z$ is a common leaf of $U$ and $T_{i+1}$ is
treated similarly.

\vskip5pt

Case 2: $u_{i+1}$ is a leaf of $T_i$. Then, by condition (II), 
$u_{i+1}$ is adjacent in
$C_i$ to a vertex $y$ such that either $d_T(u_{i+1},y)\leq p-2$, or
$d_T(u_{i+1}) = p-1$ and $y$ is not a leaf of $U$. Let $C_{i+1}$ be
the cycle obtained by replacing the edge $(y,u_{i+1})$ with the path
$yw_1R_iw_2u_{i+1}$ in $C_i$. Then $C_{i+1}$ is a Hamiltonian cycle in
$T_{i+1}^p$ because $d_T(y,w_1) = d_T(y,u_{i+1}) + 1 \leq p$ and
$d_T(w_2,u_{i+1})\leq 2$.

As in the previous case, it is easy to see that $C_{i+1}$ satisfies
condition (I). To show that $C_{i+1}$ satisfies condition (II), we
need to consider the case when $y$ is a common leaf of $U$ and
$T_{i+1}$. In this case, $y$ is also a leaf of $T_i$, so either
$d_T(y,u_{i+1})\leq p-2$, or the vertex $y'\neq u_{i+1}$
adjacent to $y$ in $C_i$ has the property of condition (II). In the
former case, we have $d_T(y,w_1) \leq p-1$, and $w_1$ is not a leaf of
$U$, and in the latter case the vertex $y'$, which is adjacent to $y$
in $C_{i+1}$, has the required property.

\vskip5pt

Since $T_\ell = T$, it follows that $C_\ell$ is a Hamiltonian cycle in
$T^p$, and hence in $G^p$, containing $v_1$,\dots,$v_k$ in order. The
lemma follows.
\end{proof}

\begin{lemma}\label{tl2}
For $t\geq 1$, let $G$ be a connected graph on at least $t$
vertices. Then there exists a map $\alpha : V(G)\rightarrow
\{1,\dots,t\}$ such that, for every two vertices $x$ and $z$ and every
integer $c\in\{1,\dots,t\}$, there
exists a sequence $x = y_0$, $y_1$, \dots, $y_{\ell-1}$, $y_\ell = z$
of distinct vertices such that $d_G(y_i, y_{i-1}) \leq t$ for $1\leq i
\leq \ell$, 
$\alpha(y_i) = c$ for $1\leq i \leq \ell-1$,
and $d_G(x,y_1) \leq t-1$ (resp.,\ $d_G(y_\ell,z)\leq t-1$)
if $\alpha(x)\neq c$ (resp.,\ $\alpha(z)\neq c$).
\end{lemma}

\begin{proof} Let $T$ be a spanning tree of $G$. Since $d_G(x,x')\leq
  d_T(x,x')$ for all $x$ and $x'$, it suffices
  to prove the lemma with $G$ replaced by $T$.
We define a map $\alpha : V(T) \rightarrow \{1,\dots,t\}$
such that, for every vertex $x$, the following conditions are
satisfied:
\begin{enumerate}
\item[(a)] there is a sequence $b_1(x)$, \dots, $b_{t-1}(x)$
  of vertices such that $d_G(x,b_i(x)) \leq
  i$ for $1\leq i \leq t-1$, and the sequence
  $\alpha(b_1(x))$,\dots,$\alpha(b_{t-1}(x))$ contains all elements of
  $\{1,\dots,t\}-\{\alpha(x)\}$;

\item[(b)] if there is a vertex $y\neq x$ such that
  $\alpha(x)=\alpha(y)$, then there is a vertex $b_t(x)\neq x$ such that
  $d_G(x,b_t(x)) \leq t$ and $\alpha(b_t(x)) = \alpha(x)$.
\end{enumerate}

Let $U_t$ be a subtree of $T$ with $|U_t|=t$. Construct a sequence
$U_t \subset U_{t+1} \subset \dots \subset U_n=T$ of subtrees of $T$
as follows: for $t\leq j \leq n-1$, choose a vertex $x_{j+1}\notin U_j$
adjacent to a vertex $y_j\in U_j$, and let $U_{j+1}$ be the tree
obtained by adjoining $x_{j+1}$ to $U_j$ by means of the edge
$(x_{j+1},y_j)$.

First, we choose a different value of $\alpha(x)$ for each $x\in
U_t$. Clearly, this assignment satisfies condition (a) above: we can
take $b_i(x)$ to be the $i$-th closest vertex of $U_t$ to $x$. 

For
$t\leq j \leq n-1$, define $\alpha(x_{j+1})$ as follows. Let
$b_1(y_j)$, \dots, $b_{t-1}(y_j)$ be the sequence described in
condition (a). Set $b_1(x_{j+1}) = y_j$ and $b_{i+1}(x_{j+1}) =
b_i(y_j)$ for $1\leq i \leq t-1$. Set $\alpha(x_{j+1}) =
\alpha(b_{t-1}(y_j))$. It is easy to verify that the sequence
$b_1(x_{j+1})$, \dots, $b_t(x_{j+1})$ satisfies conditions (a) and
(b). Also, note that if, before $\alpha(x_{j+1})$ was defined, 
$x=b_{t-1}(y_{j-1})$ was the only vertex mapped by $\alpha$
to $\alpha(x)$, then we can set $b_t(x) = x_{j+1}$ to satisfy
condition (b) for the vertex $x$.

We now prove by induction that for every $t\leq j \leq n$ and every
two vertices $x$ and $z$ of $U_j$, there exists a sequence $x =
y_0,\dots, y_\ell=z$ of vertices of $U_j$ satisfying the conditions of
the lemma. This assertion is true for $j=t$ because we can set
$y_0=x$ and $y_1 = z$ since $d_T(x,z)\leq t-1$. 

Suppose the
assertion is true for some $j$. Let $x$ and $z$ be two vertices of
$U_{j+1}$, and let $c$ be a element of $\{1,\dots,n\}$. 
Without loss of generality, assume that $z\neq x_{j+1}$. If
$x\neq x_{j+1}$, then we can find the desired sequence in $U_j$. Now,
suppose that $x=x_{j+1}$. Let $x'$ be the element of the sequence
$b_1(x)$, $\dots$, $b_t(x)$ such that $\alpha(x')=c$. Note that in our
construction, $b_t(x)$ is defined since $x\notin U_t$, so $x'$ is well defined.
By the inductive hypothesis, there is a sequence $x' = y_1, \dots,
y_{\ell} = z$ in $U_j$ such that $d_T(y_i, y_{i-1}) \leq t$ for
$2\leq i \leq \ell$, $\alpha(y_i) = c$ for $2\leq i \leq \ell-1$, and
$d_T(y_{\ell-1},z) \leq t-1$ if $\alpha(z)\neq c$. Then the sequence $x =
y_0, y_1, \dots, y_\ell=z$ satisfies the conditions of the lemma
because $d_T(x,x') \leq t$, and $d_T(x,x') \leq t-1$ if $\alpha(x) \neq
\alpha(x') = c$. Therefore, the assertion is true for $t\leq j \leq
n$, and setting $j=n$ yields the lemma.
\end{proof}

We now prove the main theorem.

\begin{thm}\label{maintheorem}
For $k\geq 3$, let $G$ be a connected graph on at least $k$
vertices. Then $G^{\lfloor 3k/2\rfloor}$ is $k$-ordered, and
$G^{\lfloor 3k/2\rfloor +1}$ is $k$-ordered Hamiltonian.
\end{thm}

\begin{proof}
Put $t = \left\lfloor {3k\over 2}\right\rfloor$. Let $v_1,\dots, v_k$
be a sequence of $k$ vertices of $G$. Let $T$ be a spanning tree of
$G$, and let $U$ be the smallest subtree of $T$ containing all the
$v_i$'s. Then all leaves of $U$ are among the $v_i$'s. 

Let $\alpha : V(U) \rightarrow \{1,\dots,t\}$ be a map satisfying the
conditions of Lemma \ref{tl2} for $G=U$.

First, we construct a cycle $C$ in $U^t$ containing the vertices $v_1$,
\dots, $v_k$ in order. Put $V=\{v_1,\dots,v_k\}$, and mark each
$v_i\in V$ for which no other element $v_j\in V$ satisfies
$\alpha(v_i) = \alpha(v_j)$. Let $S(V) = \alpha(V)$, and let $u(V)$ be
the number of unmarked vertices.

\begin{lemma}\label{justalemma2}
$|S(V)|+u(V)\leq t$.
\end{lemma}

\begin{proof}[Proof of Lemma \ref{justalemma2}]
The inequality is derived in the same way as the inequality
(\ref{eq13}) in the proof of Lemma \ref{justalemma}. Values of
$\alpha$ play the role of residues modulo $t$.
\end{proof}

Let $v_{j_1}$, \dots, $v_{j_{u(V)}}$ be the unmarked elements of
$V$. Choose a sequence $r_1$, \dots, $r_{u(V)}$ of distinct elements
of $\{1,\dots,t\}-S(V)$. Such a sequence exists because $u(V) \leq t -
|S(V)|$.
Define $r(v_i)$ as follows:
$$
r(v_i) = \left\{ \begin{array}{ll} \alpha(v_i), & \mbox{if $v_i$ is marked,}\\
                    r_\ell, & \mbox{if $v_i$ is not marked and
		      $v_i=v_{j_\ell}$}.
         \end{array}
         \right.
$$
Note that $r(v_1)$,\dots,$r(v_k)$ are all distinct since the values
of $\alpha(v_i)$ are distinct for all marked vertices $v_i$.

For $1\leq i \leq k$, let $R_i$ be the path in $U^t$ between $v_i$ and
$v_{i+1}$ traversing the sequence $v_i = y_0, \dots, y_\ell = v_{i+1}$
obtained by applying Lemma \ref{tl2} to $x=v_i$, $z=v_{i+1}$, and $c = r(v_i)$
(indices taken modulo $k$).
For $i\neq j$, we have $r(v_i)\neq r(v_j)$, so the values of $\alpha$
are different for interior vertices of $R_i$ and interior vertices of
$R_j$, hence $R_i$ and $R_j$ are interior vertex disjoint. Also,
no path $R_i$ contains a vertex $v_j$ in its interior. Indeed, if
$v_i$ is marked, then $\alpha(y) = \alpha(v_i) \neq \alpha(v_j)$ for
every interior vertex $y$ of $R_i$, so $y\neq v_j$. If $v_i$ is not
marked, then $\alpha(y) = r(v_i) \notin S(V) \ni \alpha(v_j)$ for
every interior vertex $y$ of $R_i$, so $y\neq v_j$.

It follows that joining the paths $R_1$, \dots, $R_k$ yields a cycle $C$
in $U^t$ containing the vertices $v_1$, \dots, $v_k$ in order. 
Since $C$ is a cycle in $G^t$, it follows that $G^t$ is $k$-ordered.
Note that $C$ contains the leaves of $U$ because all of them are among
the $v_i$'s. 
Also, note that for a leaf $v_i$ of $U$, we have either $r(v_{i-1}) \neq
\alpha(v_i)$ or $r(v_i)\neq \alpha(v_i)$, so the vertex adjacent to
$v_i$ in one of $R_{i-1}$ and $R_i$ is at most distance $t-1$ away
from $v_i$. Applying Lemma \ref{tl1} with $p=t+1$, we conclude that
$G^{t+1}$ is $k$-ordered Hamiltonian.
\end{proof}

\begin{cor}\label{maincor}
For a connected graph $G$ on at least $k$ vertices, the inequality
$p_k(G) \leq \lfloor 3k/2\rfloor + 1$ holds.
\end{cor}

\section{4-Ordered-Hamiltonicity of $G^4$}\label{four}

The obtained upper bound on $p_k(G)$ for an arbitrary connected graph
$G$ is not tight for small $k$. 
Indeed, Corollary \ref{maincor} states that
$p_3(G) \leq 5$ and $p_4(G) \leq 7$, whereas Corollary \ref{g3} implies
$p_3(G) \leq 3$, and our next result yields $p_4(G) \leq 4$.

\begin{thm}\label{g4}
If $G$ is a connected graph on 4 or more vertices, then $G^4$ is
4-ordered Hamiltonian.
\end{thm}

As in the proof of Theorem \ref{g3c}, it suffices to show that $T^4$
is 4-ordered Hamiltonian, where $T$ is a spanning tree of $G$. The
argument will be based on the following technical lemma which sets up
an application of Lemma \ref{tl1}.

\begin{lemma}\label{t-bar}
Let $\bT$ be a tree with at most 4 leaves, 
and let $v_1$, $v_2$, $v_3$, $v_4$ be a
sequence of 4 vertices of $\bT$ such that all leaves of $\bT$
are among $v_1$, \dots, $v_4$. Then there is a Hamiltonian cycle $C$ in
$\bT^4$ containing $v_1$, \dots, $v_4$ in the given order and
having the following property: for each vertex $v_i$ that
is a leaf of $\bT$, there is a vertex $w_i$ 
adjacent to $v_i$ in $C$ such that
$d_{\bT}(v_i, w_i) \leq 2$ or else $d_{\bT}(v_i, w_i) = 3$ and
$w_i$ is not a leaf of $\bT$.
\end{lemma}

\begin{figure}
\begin{center}
\input{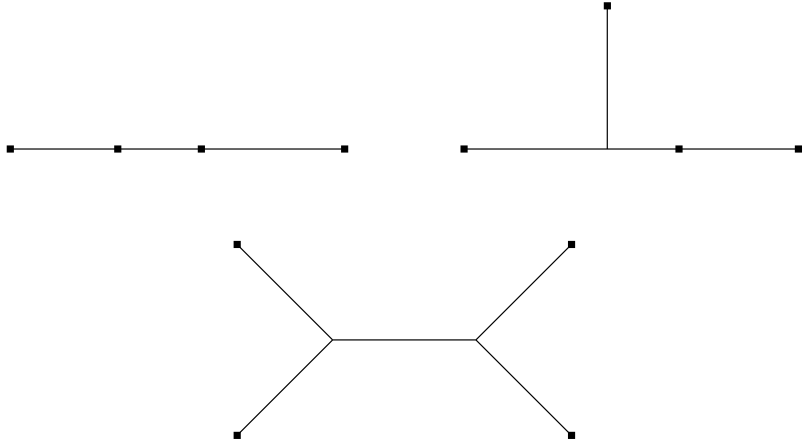}
\end{center}
\caption{The tree $\bT$ with 2, 3, and 4 leaves}
\label{types}
\end{figure}

\begin{proof}[Proof of Lemma \ref{t-bar}]
Figure \ref{types} shows the three possible shapes of $\bT$ and
locations of $v_1$, \dots, $v_4$. We will consider each case
separately.

Case 1: $\bT$ has 2 leaves, i.e. $\bT$ is a path. Let $w_1$,
$w_2$, $w_3$, $w_4$ be the vertices $v_1$, $v_2$, $v_3$, $v_4$ in the
order in which they appear in $\bT$, so that $w_1$ and $w_4$ are
the two leaves. We
will construct 6 internally disjoint paths $R_{ij}$ in $\bT^4$, $1\leq
i < j\leq 4$,
where $R_{ij}$ is a path between $w_i$ and $w_j$ and $R_{ij} = R_{ji}$
by convention.

We start by choosing integers $a$ and $b$ such that $1\leq a,
b \leq 3$ and $a + b \equiv d_{\bT}(w_2, w_3) \pmod 4$. Then we
choose an integer $c$ such that $1\leq c \leq 3$ and $c\neq
a$. Finally, we choose an integer $d$ such that $1\leq d \leq 3$,
$d\neq b$, and $c+d \not\equiv d_{\bT}(w_2, w_3) \pmod 4$.

Let $R_{14}$ be the path $w_1x_1x_2\dots x_ky_1y_2\dots y_lz_1z_2\dots
z_mw_4$ in $\bT^4$, where 
$$\begin{array}{lll}
d_{\bT}(w_1, x_1) \leq 4, &
d_{\bT}(x_i, x_{i+1}) = 4 \mbox{ for $1\leq i \leq k-1$}, &
d_{\bT}(x_k, w_2) = 4-a, \\ d_{\bT}(w_2, y_1) = a, &
d_{\bT}(y_i, y_{i+1}) = 4 \mbox{ for $1\leq i\leq l-1$}, &
d_{\bT}(y_l, w_3) = b, \\ d_{\bT}(w_3, z_1) = 4-b, &
d_{\bT}(z_i, z_{i+1}) = 4 \mbox{ for $1\leq i\leq m-1$}, & 
\mbox{and }d_{\bT}(z_m, w_4)\leq 4.
\end{array}$$

Let $R_{13}$ be the path $w_1x'_1x'_2\dots x'_py'_1y'_2\dots y'_qw_3$
in $\bT^4$,
where 
$$\begin{array}{lll}
d_{\bT}(w_1,x'_1)\leq 4, & d_{\bT}(x'_i, x'_{i+1}) = 4 \mbox{
for $1\leq i \leq p-1$}, & d_{\bT}(x'_p, w_2) = 4-c, \\
d_{\bT}(w_2, y'_1) = c, & d_{\bT}(y'_i, y'_{i+1}) = 4 \mbox{ for
$1\leq i \leq q-1$}, &
\mbox{and }d_{\bT}(y'_q, w_3) \leq 4.
\end{array}$$

Let $R_{24}$ be the path $w_2y''_1y''_2\dots y''_rz'_1z'_2\dots
z'_sw_4$ in $\bT^4$,
where 
$$\begin{array}{lll}
d_{\bT}(w_2, y''_1) \leq 4, & d_{\bT}(y''_i, y''_{i+1})
= 4 \mbox{ for $1\leq i\leq r-1$}, & d_{\bT}(y''_r, w_3) = d, \\
d_{\bT}(w_3, z'_1) = 4-d, & d_{\bT}(z'_i, z'_{i+1}) = 4 \mbox{ for
$1\leq i \leq s-1$}, & \mbox{and }d_{\bT}(z'_s, w_4) \leq 4.
\end{array}$$

Figure \ref{paths} shows the paths $R_{14}$, $R_{13}$, and
$R_{24}$. The conditions on $a$, $b$, $c$, and $d$ ensure that these
paths are internally disjoint. Notice also that we allow $k$,
$l$, $m$, $p$, $q$, $r$, and $s$ to be zero (for example, if
$d_{\bT}(w_1, w_2)\leq 4-a$, then $k=0$).

\begin{figure}
\begin{center}
\input{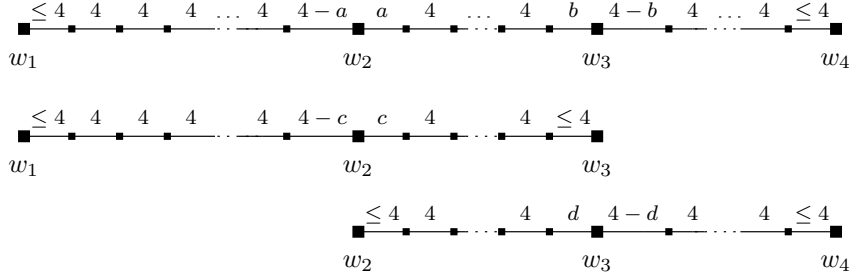}
\end{center}
\caption{The paths $R_{14}$, $R_{13}$, and $R_{24}$}
\label{paths}
\end{figure}

Let $R_{12}$, $R_{23}$, and $R_{34}$ be the paths in $\bT^4$
consisting of all vertices of $\bT$ between $w_1$ and $w_2$, $w_2$
and $w_3$, and $w_3$ and $w_4$, respectively, that are not contained
in $R_{14}$, $R_{13}$, and $R_{24}$ (the resulting paths are valid
since each of $R_{14}$, $R_{13}$, and $R_{24}$ uses every fourth
vertex on each of the three intervals). 

Now denote $v_j = w_{i_j}$ and form a cycle $C$ in $\bT^4$ 
containing $v_1$, \dots, $v_4$ in the
given order by linking the paths $R_{i_1i_2}$, $R_{i_2i_3}$,
$R_{i_3i_4}$, and $R_{i_4i_1}$.
This cycle can be easily extended to a
Hamiltonian cycle since every vertex $v\notin C$ lies between
vertices $t$ and $u$ that are adjacent in $C$, so $v$ can be inserted
between $t$ and $u$.

The conditions of the lemma are satisfied since every leaf of
$\bT$ is adjacent in $C$ to a non-leaf that is at most distance 3
away.

Case 2: $\bT$ has 3 leaves. Without loss of generality, assume that
$v_1$, $v_2$, and $v_3$ are the leaves of $\bT$. Let $v_0$ be the only
vertex of $\bT$ of degree 3. We consider two subcases.

Case 2.1: $v_4$ lies between $v_1$ and $v_0$ (the case when $v_4$ lies
between $v_3$ and $v_0$ is analogous) or $v_4 = v_0$. 
Let $P_{12}$ be the path in
$\bT$ between $v_1$ and $v_2$, and let $H$ be the Hamiltonian cycle in
$(P_{12})^2$. The part of $H$ between $v_2$ and $v_4$ that does not
contain $v_1$ contains either the vertex $v_0$ or two vertices $u$ and
$w$ that are adjacent to $v_0$ in $P_{12}$ (suppose that $u$ is closer
to $v_2$ than $w$). Let $P_{03}$ be the path
in $\bT$ between $v_0$ and $v_3$, and let $R$ be the Hamiltonian path
in $(P_{03} - v_0)^2$ between vertices $x$ and $y$, where
$d_{\bT}(v_0, x) = 1$ and $d_{\bT}(v_0, y) = 2$ (if $v_3$ is the only
vertex of $P_{03} - v_0$, then let $R$ be the ``empty'' path starting
and ending at $v_3$). Form the Hamiltonian cycle $C$ in $\bT^4$ as
follows:
start by going from $v_1$ to $v_2$ along the part of $H$ that does not
contain $v_4$, then continue along $H$ until either $u$ or the vertex
preceding $v_0$ in $H$ is
encountered, then proceed to $x$ and go to $y$ along $R$, and finally,
proceed to $w$ or $v_0$ and finish the
cycle by going to $v_1$ along the remaining part of $H$ (see Figure
\ref{cyc21}).

\begin{figure}
\begin{center}
\input{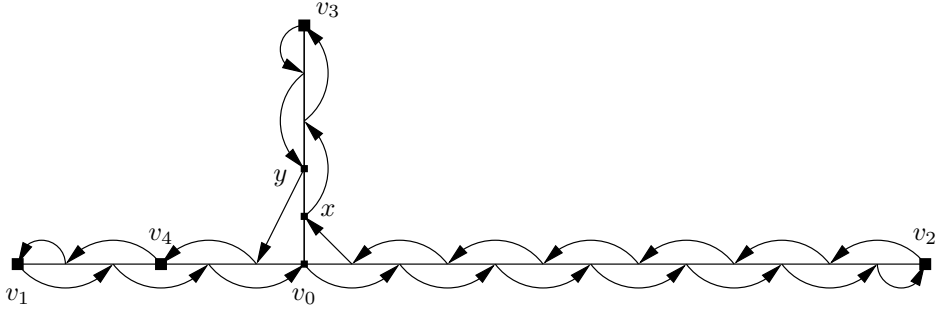}
\end{center}
\caption{The cycle $C$ in Case 2.1}
\label{cyc21}
\end{figure}

Clearly, the conditions of the lemma are satisfied since every leaf
of $\bT$ is adjacent in $C$ to a vertex of $\bT$ that is at
most distance 2 away.

Case 2.2: $v_4$ lies between $v_2$ and $v_0$. For $0\leq i,j\leq 4$
let $P_{ij}$ denote the path in $\bT$ between $v_i$ and $v_j$.   
Let $t_i$, $u_i$, and $w_i$ denote vertices of $P_{10}$, $P_{20}$, and
$P_{30}$, respectively, such that $d_{\bT}(t_i, v_0) = d_{\bT}(u_i,
v_0) = d_{\bT}(w_i,
v_0) = i$. Let $R_1$ be the Hamiltonian path in $(P_{10})^2$ between
$v_0$ and $t_1$. Let $R_2$ be the Hamiltonian path in
$(P_{42} - v_4)^2$ between $x_1$ and $x_2$, where $d_{\bT}(v_4, x_1) = 1$
and $d_{\bT}(v_4, x_2) = 2$ (if $v_2$ is the only vertex of $P_{42} -
v_4$, then let $R_2$ be the ``empty'' path starting and ending at
$v_2$). Let $R_3$ be the Hamiltonian path in $(P_{30} - v_0)^2$
between $w_1$ and $w_2$ (again, if $v_3$ is the only vertex of $P_{30}
- v_0$, then let $R_3$ be the ``empty'' path starting and ending at $v_3$).
Below is the procedure to construct a Hamiltonian cycle $C$ in $\bT^4$
with the
desired properties.

\begin{figure}
\begin{center}
\input{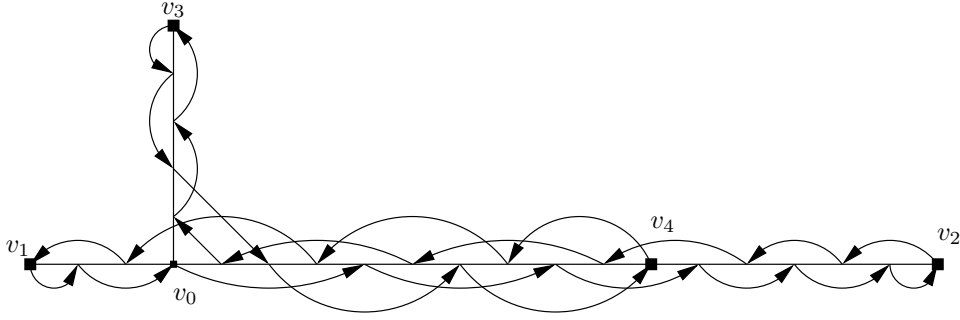}
\end{center}
\caption{The cycle $C$ in Case 2.2}
\label{cyc22}
\end{figure}

\begin{enumerate}

\item Start at $v_1$. If $d_{\bT}(v_0, v_4)$ is divisible by 4, then go
   along $R_1$ to $t_1$, otherwise go along $R_1$ to $v_0$.

\item If $d_{\bT}(v_0, v_4) \geq 4$, then use the path
   $t_1u_3u_7u_{11}\dots$ or the path $v_0u_4u_8u_{12}\dots$ to get
   within distance 3 of $v_4$.

\item Proceed to $x_1$ and 
   go along $R_2$ to $x_2$ (or stay at $x_1$ if $x_1 = v_2$).

\item If $d_{\bT}(v_0, v_4) \geq 2$, 
   then proceed to either $u_m$ or $u_{m-1}$, where
   $d_{\bT}(v_4, u_m) = 1$ and $d_{\bT}(v_4, u_{m-1}) = 2$ (one of
   $u_m$ and $u_{m-1}$ has not been used yet) and use the path
   $u_mu_{m-4}u_{m-8}\dots$ or $u_{m-1}u_{m-5}u_{m-9}\dots$ to get to
   one of $u_1$, $u_2$, or $u_3$ (notice that it is impossible to
   encounter $v_0$ at this step since $v_0$ was encountered in Step 1
   unless $d_{\bT}(v_0, v_4)$ is divisible by 4, in which case $u_2$
   or $u_3$ will be encountered).

\item Proceed to $w_1$ and 
   go along $R_3$ to $w_2$ (or stay at $w_1$ if $w_1 = v_3$).

\item If $d_{\bT}(v_0, v_4) \geq 4$, then
   proceed to $v_0$, $u_1$, or $u_2$ (one of these has not been used
   yet), and follow the path $v_0u_4u_8\dots$, $u_1u_5u_9\dots$, or
   $u_2u_6u_{10}\dots$ to get within distance 4 of $v_4$.

\item Proceed to $v_4$. If $d_{\bT}(v_0, v_4) \leq 3$, then go to Step 9.

\item One of the vertices $u_m$, $u_{m-1}$, $u_{m-2}$, $u_{m-3}$ has not
   been used yet; let $u_{m-j}$ be this vertex ($0\leq j \leq 3$). Use
   the path $v_4u_{m-j}u_{m-j-4}u_{m-j-8}\dots$ to get to one of
   $v_0$, $u_1$, $u_2$, $u_3$.

\item If $v_0$ was encountered in Step 1, then proceed to $t_1$ and finish
   the cycle by going along the remaining part of $R_1$ to
   $v_1$. Otherwise, proceed to $v_0$ and go along the remaining part
   of $R_1$ to $v_1$.

\end{enumerate}

Let us verify that the conditions of the lemma are satisfied. The leaf
$v_1$ is adjacent in $C$ to two vertices of the path $P_{12}$, at least one
of which is a non-leaf that is at most distance 3 away. The same is
true of $v_2$ unless $d_{\bT}(v_2,v_4) \leq 2$ 
and $d_{\bT}(v_0, v_4) = 1$, in which
case $v_2$ is adjacent in $C$ to $v_0$ or $x_1$, one of 
which is at most distance 2
away from $v_2$. Finally, $v_3$ is adjacent in $C$ to two vertices of the path
$P_{32}$, one of which is a non-leaf that is at most distance 3 away.

Case 3: $\bT$ has 4 leaves $v_1$, $v_2$, $v_3$, and $v_4$. There is
either one vertex of degree 4, which we denote by $u_0$, or two
vertices of degree 3, which we denote by $u_0$ and $u_m$, where $m =
d_{\bT}(u_0, u_m)$. Let $P = u_0u_1u_2\dots u_m$ be the path in $\bT$
between $u_0$ and $u_m$ (if the degree of $u_0$ is 4, then let $P$ be
the ``empty'' path consisting of $u_0$ alone). Without loss of
generality, 
assume that
$u_0$ is closer to $v_1$ than $u_m$. For $1\leq i \leq 4$, let $P_i$
be the path in $\bT$ between $v_i$ and
the closest endpoint of $P$. Let $x^1_i$, $x^2_i$, $x^3_i$, and $x^4_i$ denote 
the vertex of $P_1$, $P_2$, $P_3$, and $P_4$, respectively, the
distance from which to the closest endpoint of $P$ is $i$.

Case 3.1: $v_1$ and $v_2$ are on one side of $P$, and $v_3$ and $v_4$
are on the other side of $P$ (the case when $v_1$ and $v_4$ are on one
side of $P$ is analogous). 
For $1\leq i \leq 2$, let $R_i$ be the Hamiltonian path in $(P_i)^2$
between $x^i_1$ and $u_0$, and for $3\leq i \leq 4$ let $R_i$ be the
Hamiltonian path in $(P_i-u_m)^2$ between $x^i_1$ and $x^i_2$ (if
$x^i_1 = v_i$, then let $R_i$ be the ``empty'' path consisting of $v_i$
alone). Below
is the procedure to construct a Hamiltonian cycle $C$ 
in $\bT^4$ with the desired properties.

\begin{figure}
\begin{center}
\input{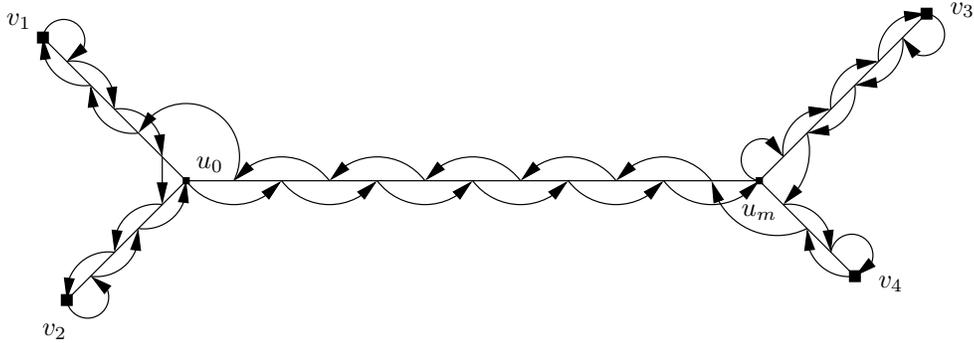}
\end{center}
\caption{The cycle $C$ in Case 3.1}
\label{cyc31}
\end{figure}

\begin{enumerate}

\item Start at $v_1$ and go along $R_1$ to $x^1_1$.

\item Proceed to $x^2_1$ and go along $R_2$ to $u_0$.

\item If $m \geq 2$, then use the path $u_0u_2u_4\dots$ to get within
   distance 1 of $u_m$.

\item Proceed to $x^3_1$ and go along $R_3$ to $x^3_2$ (or stay at
   $x^3_1$ if $x^3_1 = v_3$).

\item Proceed to $x^4_1$ and go along $R_4$ to $x^4_2$ (or stay at
   $x^4_1$ if $x^4_1 = v_4$).

\item If $m\geq 1$, then proceed to $u_m$ or $u_{m-1}$ (one of these
   vertices has not been used yet) and use the path $u_m u_{m-2}
   u_{m-4}\dots$ or the path $u_{m-1} u_{m-3} u_{m-5}\dots$ to get
   within distance 1 of $u_0$.

\item Proceed to $x^1_1$ and go along the remaining part of $R_1$ to
   $v_1$.

\end{enumerate}

Let us verify that the conditions of the lemma are satisfied. If
the distance between $v_i$ and the closest endpoint of $P$ is greater
than 1, then $v_i$ is adjacent in $R_i$, and hence in $C$, to a vertex
that is at most distance 2 away. If $d_{\bT}(v_1, u_0) = 1$, then $v_1$
is adjacent in $C$ to $x^2_1$, which is distance 2 away. If
$d_{\bT}(v_2, u_0) = 1$, then $v_2$ is adjacent in $C$ to $u_0$. If
$d_{\bT}(v_3, u_m) = 1$, then $v_3$ is adjacent in $C$ to $u_m$ or
$u_{m-1}$, which are distance 1 and 2 away, respectively. If
$d_{\bT}(v_4, u_m) = 1$ and $m \geq 1$, then $v_4$ is adjacent in $C$ to either
$u_m$ or $u_{m-1}$, which are distance 1 and 2 away,
respectively. Finally, if $d_{\bT}(v_4, u_m) = 1$ and $m=0$, then $v_4$
is adjacent in $C$ to $x^1_1$, which is distance 2 away.

Case 3.2: $v_1$ and $v_3$ are on one side of $P$, and $v_2$ and $v_4$
are on the other side of $P$. 
%
%
%

Let $R_1$ be the Hamiltonian path in $(P_1)^2$ between $x^1_1$ and
$u_0$. For $2\leq i \leq 4$, let $R_i$ be the Hamiltonian path in $(P_i
- u_j)^2$ between $x^i_1$ and $x^i_2$, where $u_j$ is the endpoint of 
$P_i$ different from $v_i$ (if $v_i = x^i_1$, then let $R_i$ be the
``empty'' path consisting of $v_i$ alone).
Below is the procedure to construct a Hamiltonian cycle $C$ 
in $\bT^4$ with the desired properties.

\begin{figure}
\begin{center}
\input{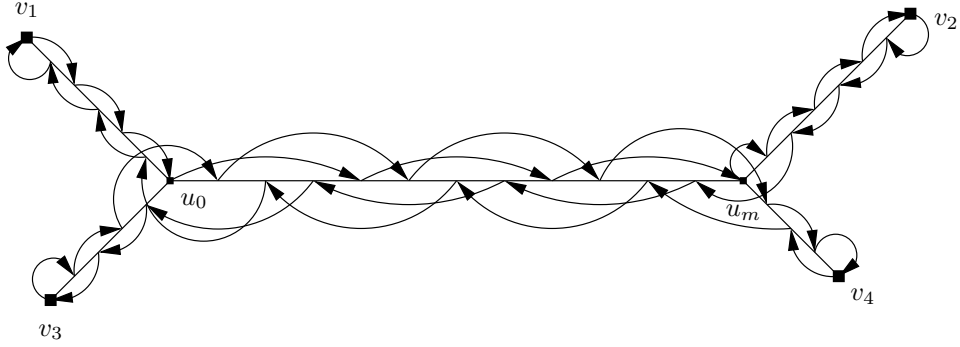}
\end{center}
\caption{The cycle $C$ in Case 3.2}
\label{cyc32}
\end{figure}

\begin{enumerate}

\item Start at $v_1$ and go along $R_1$ to $u_0$.

\item If $m \geq 4$, then use the path $u_0u_4u_8\dots$ to get within
   distance 3 of $u_m$.

\item Proceed to $x^2_1$ and go along $R_2$ to $x^2_2$ (or stay at
   $x^2_1$ if $x^2_1 = v_2$).

\item If $m \geq 2$, then proceed to $u_m$ or $u_{m-1}$ (one of them has
   not been used yet) and use the path $u_m u_{m-4} u_{m-8}\dots$ or
   the path $u_{m-1}u_{m-5}u_{m-9}\dots$ to get within distance 3 of
   $u_0$.

\item Proceed to $x^3_1$ and go along $R_3$ to $x^3_2$ (or stay at
   $x^3_1$ if $x^3_1 = v_3$).

\item If $m \geq 2$, then proceed to $u_1$ or $u_2$ (one of them
   has not been used yet) and use the path $u_1 u_5 u_9
   \dots$ or the path $u_2 u_6 u_{10}\dots$ to get within distance 3 of $u_m$.
   Otherwise, if $m=1$, then proceed to $u_1$ (in this case we proceeded from
   $x^2_2$ or $x^2_1$ directly to $x^3_1$, so $u_1$ has not been used).

\item If $m\geq 3$ and $u_{m-3}$ is encountered at the end of Step 6 or
   if $m=0$, then proceed to
   $x^4_1$ and go along $R_4$ to $x^4_2$ (or stay at $x^4_1$ if $x^4_1
   = v_4$), otherwise proceed
   to $x^4_2$ and go along $R_4$ to $x^4_1$ (or proceed to $x^4_1$ and
   stay there if $x^4_1 = v_4$).

\item If $m\geq 3$, then 
   one of the vertices $u_m$, $u_{m-1}$, $u_{m-2}$, $u_{m-3}$ has not been
   encountered yet; denote this vertex by $u_{m-j}$ (if the last
   encountered vertex is $x^4_2$, then $u_{m-3}$ has been already
   encountered, so $j\leq 2$). Proceed to $u_{m-j}$ and use the path
   $u_{m-j}u_{m-j-4}u_{m-j-8}\dots$ to get within distance 3 of $u_0$.

\item Proceed to $x^1_1$. (If $1 \leq m \leq 2$, then the last encountered
    vertex is $x^4_1$, so $x^1_1$ is at most distance 4 away, and if
    $m=0$, then $d_{\bT}(x^4_2, x^1_1) = 3$.)

\item Finish the cycle by going along the remaining part of $R_1$ to $v_1$.

\end{enumerate}

Let us verify that the obtained cycle $C$ satisfies the conditions of
the lemma. Let $P_{ij}$ be the path in $\bT$ between $v_i$ and
$v_j$. If $m\geq 1$, then $v_1$ and $v_4$ are each adjacent in $C$ to
two vertices of $P_{14}$, one of which is a non-leaf that is
at most distance 3 away, and $v_2$ and $v_3$ are each adjacent to two vertices
of $P_{23}$, one of which is a non-leaf that is at most distance 3
away. Now suppose that $m = 0$. If $d_{\bT}(v_i, u_0) > 1$, then $v_i$
is adjacent in $R_i$, and hence in $C$, to a vertex of $\bT$ that is
at most distance 2 away. If $d_{\bT}(v_1, u_0) = 1$, then $v_1$ is adjacent in
$C$ to $u_0$. If $d_{\bT}(v_2, u_0) = 1$, then $v_2$ is adjacent in $C$
to $x^3_1$, which is distance 2 away. If $d_{\bT}(v_3, u_0) = 1$, then  
$v_3$ is adjacent in $C$ to $x^4_1$, which is distance 2
away. Finally, if $d_{\bT}(v_4, u_0) = 1$, then $v_4$ is adjacent in
$C$ to $x^1_1$, which is distance 2 away.
\end{proof}

Theorem \ref{g4} now follows from Lemma \ref{t-bar} and Lemma
\ref{tl1}.

\section{Powers of cycles}\label{cycles}

In this section, we compute $p_5(C_n)$ and give a lower bound on $p_k(C_n)$.

\begin{thm}\label{cyc3}
Let $C_n$ denote the cycle on $n$ vertices. If $n\geq 5$, then
$(C_n)^3$ is 5-ordered Hamiltonian.
\end{thm}

\begin{proof}
Let $v_1$, $v_2$, \dots, $v_5$ be a sequence of 5 vertices of $C_n$,
and let $w_1$, $w_2$, \dots, $w_5$ be the same 5 vertices in the order
in which they appear in the cycle $C_n$. For $1\leq i \leq 5$, let
$P_i$ denote the portion of $C_n$ between $w_i$ and $w_{i+1}$
containing no other vertices of the sequence $w_1$, \dots, $w_5$
(indices taken modulo 5). We will construct 10 internally-disjoint
paths $R_{ij}$ in $(C_n)^3$, where $R_{ij}$ is a path between $w_i$
and $w_j$ for $1\leq i<j \leq 5$. We adopt the convention that $R_{ji}
= R_{ij}$.

\begin{figure}
\begin{center}
\input{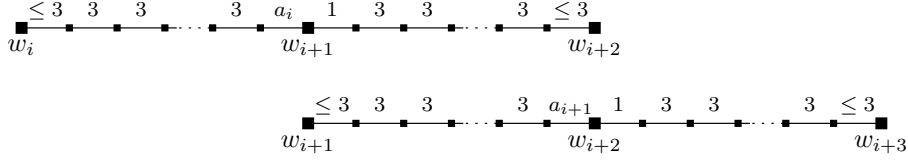}
\end{center}
\caption{The paths $R_{i,i+2}$ and $R_{i+1,i+3}$}
\label{paths_c}
\end{figure}

First, we construct $R_{i,i+2}$ for $1\leq i \leq 5$. Let $a_1$,
$a_2$, \dots, $a_5$ be integers such that $a_i \in \{1,2\}$ and $a_i +
1 \not\equiv |P_i| - 1 \pmod 3$ for all $i$
(note that $|P_i| - 1$ is the length of $P_i$).

For $1\leq i \leq 5$, let $R_{i,i+2}$ be the path $w_it^i_1t^i_2\dots
t^i_{l_i} u^i_1 u^i_2 \dots u^i_{m_i} w_{i+2}$, where $t^i_j \in P_i$,
$u^i_j \in P_{i+1}$, and
$$
\begin{array}{ll}
d_{P_i}(w_i, t^i_1) \leq 3, & d_{P_i}(t^i_j, t^i_{j+1}) = 3 \mbox{ for
  $1\leq j \leq l_i - 1$}, \\ d_{P_i}(t^i_{l_i}, w_{i+1}) = a_i, &
d_{P_{i+1}}(w_{i+1}, u^1_1) = 1, \\ d_{P_{i+1}}(u^i_j, u^i_{j+1}) =
  3 
\mbox{ for
  $1\leq j \leq m_i - 1$}, & \mbox{and }d_{P_{i+1}}(u^i_{m_i}, w_{i+2}) \leq 3.
\end{array}
$$
We allow $l_i = 0$ if $d_{P_i}(w_i,w_{i+1})\leq a_i$ and $m_i = 0$ if
$d_{P_{i+1}}(w_{i+1}, w_{i+2}) = 1$. Figure \ref{paths_c} shows the
paths $R_{i,i+2}$ and $R_{i+1,i+3}$. Since $a_{i+1} + 1\not\equiv
|P_{i+1}| - 1 \pmod 3$, these paths are internally disjoint.

For $1\leq i \leq 5$, let $R_{i,i+1}$ be the path in $(C_n)^3$
between $w_i$ and $w_{i+1}$ containing all vertices of $P_i$ that are
not contained in $R_{i, i+2}$ and $R_{i-1,i+1}$.

Define $i_j$ by $v_j = w_{i_j}$ and form a cycle $C$ in $(C_n)^3$
containing $v_1$, \dots, $v_5$ in order by linking the paths
$R_{i_1i_2}$, $R_{i_2i_3}$, \dots, $R_{i_5i_1}$. Finally, we extend
$C$ to a Hamiltonian cycle using the following procedure. 

\begin{enumerate}

\item If $C$ is
not Hamiltonian, choose a vertex $z\notin C$ of $C_n$ adjacent in $C_n$ to a
vertex $u\in C$.

\item At most one vertex adjacent to $u$ in
   $(C_n)^3$ is more than distance 3 away from $z$. Since $u$ is
   adjacent to two vertices in $C$, 
   we can choose a vertex $t$ adjacent to $u$ in $C$ such
   that $d_{C_n}(t,z) \leq 3$.

\item Replace the edge $(t,u)$ of $C$ with the path $tzu$.

\item If $C$ is not Hamiltonian, return to Step 1.

\end{enumerate}

Since during the procedure we insert the remaining vertices
into $C$ without changing the order of the vertices already in $C$,
the order in which $v_1$, \dots, $v_5$ are contained in $C$ is
preserved, so $(C_n)^3$ is 5-ordered Hamiltonian. This completes the proof of
the theorem.
\end{proof}

\begin{cor}\label{ham3}
If $G$ is a Hamiltonian graph on 5 or more vertices, then $G^3$ is
5-ordered Hamiltonian.
\end{cor}

\begin{proof}
Let $C$ be a Hamiltonian cycle in $G$. Then $C^3$ is 5-ordered
Hamiltonian. Since $C^3$ is a subgraph of $G^3$ with the same vertex
set, it follows that $G^3$ is 5-ordered Hamiltonian, too.
\end{proof}

The next proposition yields a lower bound on $p_k(C_n)$.

\begin{prop}\label{sharp2}
If $m\geq 3$ and
$n$ is sufficiently large,
then $(C_n)^m$ is not $2m$-ordered.
\end{prop}

\begin{figure}
\begin{center}
\input{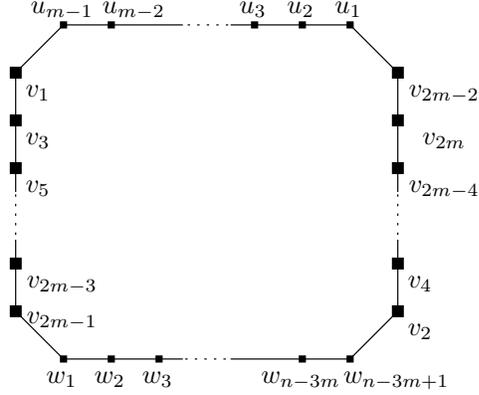}
\end{center}
\caption{The cycle $C_n$}
\label{cm}
\end{figure}

\begin{proof}
Let $v_1$, $v_3$, $v_5$, \dots, $v_{2m-1}$, 
$w_1$, $w_2$, \dots, $w_{n-3m+1}$, 
$v_2$, $v_4$, \dots, $v_{2m-4}$, $v_{2m}$, $v_{2m-2}$,
$u_1$, $u_2$, \dots,
$u_{m-1}$
be the vertices of $C_n$ (in this order; see Figure
\ref{cm}). 
Suppose
$(C_n)^m$ is $2m$-ordered; then there
is a cycle $C$ in $(C_n)^m$ that contains $v_1$, $v_2$, \dots,
$v_{2m}$ in order. For $1\leq i \leq 2m$, let $R_i$ be the portion of
$C$ between $v_i$ and $v_{i+1}$ that does not contain any other
vertices of the sequence $v_1$, \dots, $v_{2m}$ (indices taken modulo $2m$).
Denote $U = \{u_1, u_2, \dots, u_{m-1}\}$ and $W = \{w_1, w_2, \dots,
w_{n-3m+1}\}$.

If $n$ is sufficiently large, then we have $d_{C_n}(v_i, v_{i+1}) > m$
for all $i$. Therefore, each $R_i$ must contain a vertex of $U$ or a
vertex of $W$. Since $|U| = m-1$, there are at least $m+1$ paths $R_i$
that do not contain vertices of $U$ and hence consist only of vertices
in $\{v_i, v_{i+1}\} \cup W$. Each such $R_i$ must contain
at least $\left\lfloor{|W|\over m}\right\rfloor$ 
vertices of $W$, otherwise there would be $m$ consecutive vertices of
$W$ not contained in $R_i$, and hence
$R_i$ would have two adjacent vertices that are more than distance $m$ apart
in $C_n$. We conclude that at least $m+1$ of the paths $R_i$ have at
least $\left\lfloor{|W|\over m}\right\rfloor$ interior vertices, thus
$$
|C_n| > (m+1) \left\lfloor{|W|\over m}\right\rfloor = (m+1)
\left\lfloor{n-3m+1\over m}\right\rfloor > n
$$
if $n$ is sufficiently large --- a contradiction. It follows that
$(C_n)^m$ is not $2m$-ordered.
\end{proof}

\begin{cor}
$p_k(C_n) \geq \lfloor k/2\rfloor+1$ for sufficiently large $n$.
\end{cor}

\section{Conclusion}\label{conclusion}

We conclude the paper by giving a brief summary of the results and
posing questions for further research.
We have shown that $p_k(G)\leq \lfloor 3k/2 \rfloor+1$
for any connected graph $G$ on at least $k$ vertices. This bound is
almost tight because $p_k(P_n) \geq \lfloor 3k/2 \rfloor - 2$ for
sufficiently large $n$. 

\begin{q}
What is the best possible upper bound on $p_k(G)$ for an arbitrary
connected graph?
\end{q}

For $k=4$, the optimal bound  $p_4(G) \leq 4$ follows from Theorem
\ref{g4}, and tightness follows from Theorem \ref{pathlowerbound}.

In the case of the path $P_n$, the value of $p_k(P_n)$ obeys the inequalities
$\lfloor 3k/2\rfloor -2 \leq p_k(P_n) \leq \lfloor 3k/2 \rfloor -
1$ for sufficiently large $n$. 
Thus $p_k(P_n)$ can take one of two possible values. Note that
$p_3(P_n) = 3$ and $p_4(P_n)=4$, so both inequalities can be tight for
small $k$.

\begin{q}
What is the exact value of $p_k(P_n)$?
\end{q}

Another interesting problem is to advance our results on the value of
$p_k(C_n)$. Theorem \ref{cyc3} and Corollary
\ref{sharp2} imply $p_5(C_n)=3$.

\begin{q}
Determine the exact value or give bounds on $p_k(C_n)$.
\end{q}

Finally, a more difficult task is to generalize Theorem
\ref{g2} to imply $k$-ordered Hamiltonicity of powers of
$2$-connected graphs.

\begin{q}
Give bounds on $p_k(G)$ for a $2$-connected graph $G$ on at least $k$ vertices.
\end{q}

\section{Acknowledgments}\label{acknowledgments}

This research was begun the 2002 Research Experience for
Undergraduates at the University of Minnesota Duluth. I
would like to thank the program director Joseph Gallian for suggesting
this problem and encouraging me to work on it. I would also like to
thank advisors Geir Helleloid and Phil Matchett for many useful
comments on this paper. Finally, I want to thank all other student
participants and visitors for a wonderful summer.

\end{document}